# The realization of positive random variables via absolutely continuous transformations of measure on Wiener space[*]

**D. Feyel**

*Université d'Evry-Val-d'Essone, 91025 Evry Cedex, France*
*e-mail:* `Denis.Feyel@univ-evry.fr`

**A. S. Üstünel**

*ENST, Dépt. Infres, 46 rue Barrault, 75634, Paris Cedex 13, France*
*e-mail:* `ustunel@enst.fr`

**M. Zakai**

*Department of Electrical Engineering*
*Technion–Israel Institute of Technology, Haifa 32000, Israel*
*e-mail:* `zakai@ee.technion.ac.il`

**Abstract:** Let $\mu$ be a Gaussian measure on some measurable space $\{W = \{w\}, \mathcal{B}(W)\}$ and let $\nu$ be a measure on the same space which is absolutely continuous with respect to $\nu$. The paper surveys results on the problem of constructing a transformation $T$ on the $W$ space such that $Tw = w + u(w)$ where $u$ takes values in the Cameron-Martin space and the image of $\mu$ under $T$ is $\mu$. In addition we ask for the existence of transformations $T$ belonging to some particular classes.

## Contents

---

[*]This is an original survey paper

## 1. Introduction

The main topic addressed in this paper can be outlined as follows: Let $(W, \mathcal{B}(W))$ be a measurable space on which we are given two probability measures $\rho$ and $\nu$. Assume that $\nu$ is absolutely continuous with respect to $\rho$. Is there a measurable transformation $T : W \to W$ such that the image of $\rho$ under $T$ is $\nu$? Of course at this generality the answer will be too vague to be of some interest, hence let us provide the space $(W, \mathcal{B}(W))$ with more analytical structure in such a way that we can construct at least a Gaussian measure on it.

We shall call the triple $(W, H, \mu)$ an abstract Wiener space if $W$ is a separable Fréchet space, $H$ is a separable Hilbert space, called the Cameron-Martin space, which is densely and continuously embedded in $W$. We also identify $H$ with its continuous dual, hence the imbeddings $W^* \hookrightarrow H^* = H \hookrightarrow W$ are all dense and continuous. $(W, \mathcal{B}(W))$ is equipped with a probability measure $\mu$ such that for every $e$ in $W^*$, ${}_{W^*}\langle e, w \rangle_W$ has the law $N(0, |e|_H^2)$ under $\mu$. Although in the literature $W$ is in general supposed to be a Banach space, it is more interesting to assume $W$ to be a (separable) Fréchet space, since such a space is a countable projective limit of Banach spaces, all the properties of abstract Banach-Wiener spaces extend immediately to our case. Another justification of this choice lies in the fact that neither the classical Wiener space $C_0(\mathbb{R}_+, \mathbb{R})$ nor $\mathbb{R}^{\mathbb{N}}$ are Banach spaces (of course they are Fréchet!) although we need them badly in stochastic analysis[1]. Let $T$ be a measurable transformation of $W$. The problem of determining whether the induced measure $\nu = T\mu$, where $T\mu(A) = \mu(T^{-1}A)$, is absolutely continuous with respect to $\mu$ and to find the corresponding Radon-Nikodym derivative, in particular the case where $T$ is a perturbation of identity, i.e., $w \to T(w) = w + u(w)$ where $u : W \to H$, has been treated extensively (cf. [23] and the references therein). In this paper we consider the converse problem; namely, given a positive random variable $L$ on $(W, H, \mu)$ such that and $E[L] = 1$, is there a measurable transformation $T$ on $W$ such that $(dT\mu/d\mu) = L$? In addition to asking for *any* $T$, we will also ask for the existence of transformation $T$ belonging to some particular classes which answer to the latter question. To be more specific we set:

**Definition 1.1.** *Let $(W, H, \mu)$ be an abstract Wiener space and let $\nu$ be a measure on $W$ which is absolutely continuous with respect to $\mu$, i.e. $d\nu = L \, d\mu$, (or a r.v. $L$ such that $L \geq 0$ a.s. and $E_\mu[L] = 1$), then:*

---

[1] Readers who prefer to consider the classical Wiener space may, throughout the paper, suppose $W$ to be the standard $d$-dimensional Wiener space $C_0([0, 1], \mathbb{R}^d)$ with (a) $(t, w) \to w(t) = W_t(w)$ being the Wiener process and $W$ endowed with the supremum norm $\|w\| = \sup_{t \in [0,1]} |w(t)|_{\mathbb{R}^d}$ and (b) The Cameron Martin space being the space of $\mathbb{R}^d$ valued absolutely continuous functions $h(t)$ where $h(t) = \int_0^t \dot{h}(t) dt$ with the norm $|h|_H^2 = \int_0^1 |\dot{h}(t)|^2 \, dt$.

a) The measure $\nu$ (or the r.v. $L$) will be said to be realizable or representable if there exists a measurable transformation $T$ such that $T\mu = \nu$. Such a $T$ will be called a realization of $L = d\nu/d\mu$. This is equivalent to the validity of the following:
$$E_\mu[g \circ T] = E_\mu[g\,L]\,,$$
for any nice $g$ on $W$.

b) The measure $\nu \ll \mu$ or the corresponding Radon-Nikodym derivative will be said to be "realizable (or representable) by a shift" if there exists an $H$-valued random variable $u$ (a "shift") such that for $Tw = w + u(w)$, $T\mu = \nu$.

In the class (b) we will also consider several particular types of realizations:

(b1) The shift $u$ is adapted to some particular filtration on $W$.
(b2) The transformation $T$ is "triangular" (to be defined in section 5).
(b3) The shift $u$ is the gradient of some random variable $\varphi$.

In addition to these realization problems we will also consider two closely related realization problems which will be described at the end of this section.

We turn now to an outline of the paper. In the next section we summarize the notions of gradient and divergence on the Wiener space and some associated relations. In section 3 we consider the problem of the realization by adapted shifts, i.e. case of (b-1). This was introduced and treated by Üstünel and Zakai cf. in Section 2.7 of [23]. Consider first the classical setup of the Wiener space, let $\mathcal{F}_t$ to be the subsigma field induced by $\{w_\theta, \theta \leq t\}$ and consider the realization of $L$ by $Tw = w + u(w)$ where $u$ is required to be a.s. in $H$ and adapted to the filtration $(\mathcal{F}_t, t \in [0,1])$. In the abstract Wiener space there is no natural filtration, however, we can adjoin to the AWS a filtration induced by a resolution of the identity on $H$ (cf. e.g. [23]). In [23], among other results, it is proven that given $L$, we can find a sequence of adapted $H$-valued shifts, say $u_n$, such that measures induced by them, and the corresponding Radon-Nikodym derivatives of $T_n\mu$ (w.r. to $\mu$) $L_n$, satisfy

$$\lim_{n\to\infty} E[|L - L_n|] = 0$$

or

$$\sup \sum_i \Big|\nu(A_i) - T_n\mu(A_i)\Big| \to 0$$

and the supremum is taken over all finite partitions of $W$.

In section 4 we present the result of Fernique [7]; in our context it states that for any $\nu \ll \mu$ then there is a probability space $(\Omega, \mathcal{F}, P)$ and three random variables $T, B$ and $Z$ defined on this space such that $T$ and $B$ are $W$-valued, $Z$ is $H$-valued, $T(P) = \nu$, $B(P) = \mu$ and $T = B + Z$.

Section 5 represents the work of Bogachev, Koleshnikov and Medvedev [1, 2] and is based on a transformation of $\mu$ to $\nu$ by a "triangular transformation". This approach yields not only the result that every measure that is absolutely continuous with respect to the abstract Wiener space is realizable, it yields also that it is realizable by a shift *without* any further restrictions on $L$.

Section 6 follows the results of Feyel and Üstünel [9, 10]. It has its origin in the modern theory of optimal transportation (cf. e.g. [3, 11, 17, 18, 24]) i.e., given an absolutely continuous probability (w.r. to the Lebesgue measure) on $\mathbb{R}^n$, it can be transformed into any probability measure $\beta$ on $\mathbb{R}^n$ by a transformation $T(x) = \nabla \Phi(x)$ where $\Phi$ is a convex function, and this transformation is optimal in a certain sense. In [9] and [10] Feyel and Üstünel extended these transportation results to the infinite dimensional setting and some of these results are briefly explained here. In particular, the case of Wiener measure can be formulated as follows: if a probability measure $\nu$ on $W$ is at finite Wasserstein distance from the Wiener measure $\mu$, where the Wasserstein distance is defined with respect to the Hilbert norm of the Cameron-Martin space $H$ (which is not continuous on $W$), then there exists a Wiener function $\varphi$ which is Sobolev differentiable such that the transformation defined by

$$T(w) = w + \nabla \varphi(w)$$

maps $\mu$ to $\nu$ (note that there is no necessity for the absolute continuity of $\nu$ w.r. to $\mu$). We indicate some applications of these results to the polar factorization of the absolutely continuous perturbations of identity. We further extend these results by defining the notion of measures which are at locally finite Wasserstein distance from each other and prove that if at least one of them is a *spread* measure, then the other one can be written as the image of the spread one under a locally cyclically monotone map of the form of a perturbation of identity.

Section 7 deals with the following problem. Let $y = u(w') + w$ where $w$ and $w'$ represent *independent* Wiener processes where $u(w')$ takes values in the Cameron-Martin space, i.e. $u$ represent a "signal" and $w$ is an additive "noise" which is independent of the signal. Let $\nu_y$ denote the measure induced by $y$ on $W$ and let $L(w)$ denote the Radon-Nikodym derivative (or likelihood ratio) $d\nu_y/d\mu$. This model is relevant to many communication and information theory problems which makes it interesting to consider the converse problem: given a positive random variable $L(w)$ on $(W, H, \mu)$ satisfying $E[L] = 1$, under what conditions does there exist a representation of "signal plus independent noise" realization such that the Radon-Nikodym derivative of $\nu_y$ with respect to $\mu$ is $L$. This problem was first considered in [25] and further clarified in [16].

**Remark:** A problem which is also related to the problems considered in this paper is that of the realization of the random variable $L(w) = 1$, i.e. the class of measure preserving transformation on the Wiener space. This problem is treated in [23] and [26] and will not be discussed in this paper.

## 2. Preliminaries

*The gradient:* Let $(W, H, \mu)$ be an AWS and let $e_i, i = 1, 2, \ldots$ be a sequence of elements in $W^*$. Assume that the image of $e_i$ in $H$ form a complete orthonormal base in $H$. Let $f(x_1, \ldots, x_n)$ be a smooth function on $\mathbb{R}^n$ and denote by $f'_i$ the partial derivative of $f$ with respect to the $i$-th coordinate. We will denote $_{W^*}\langle e, w\rangle_W$ by $\delta e$.

For cylindrical smooth random variables $F(w) = f(\delta e_1, \ldots, \delta e_n)$, define $\nabla_h F = \frac{dF(w+\varepsilon h)}{d\varepsilon}\Big|_{\varepsilon=0}$. Therefore we set the following: $\nabla_h F = (\nabla F, h)$ where $\nabla F$, the gradient, is $H$-valued. For $F(w) = \delta e$, $\nabla F = e$, and

$$\nabla F = \sum_{i=1}^n f_i'(\delta e_1, \ldots, \delta e_n) \cdot e_i. \tag{2.1}$$

It can be shown that this definition is closable in $L^p(\mu)$ for any $p > 1$. Hence it can be extended to a wider class of functions. We will restrict ourselves to $p = 2$, consequently the domain of the $\nabla$ operation can be extended to all functions $F(w)$ for which there exists a sequence of smooth cylindrical functions $F_m$ such that $F_m \to F$ in $L_2$ and $\nabla F_m$ is Cauchy in $L_2(\mu, H)$. In this case set $\nabla F$ to be the $L_2(\mu, H)$ limit of $\nabla F_m$. This class of r.v. will be denoted $\mathbb{D}_{2,1}$. It is a closed linear space under the norm

$$\|F\|_{2,1} = E_\mu^{\frac{1}{2}}(F)^2 + E_\mu^{\frac{1}{2}}|\nabla F|_H^2. \tag{2.2}$$

Similarly let $K$ be an Hilbert space and $k_1, k_2, \ldots$ a complete orthonormal base in $K$. Let $\varphi$ be the smooth $K$-valued function $\varphi = \sum_{j=1}^m f_j(\delta e_1, \ldots, \delta e_n) k_j$ define

$$\nabla \varphi = \sum_{j=1}^m \sum_{i=1}^n (f_j)_i'(\delta e_1, \ldots, \delta e_n) e_i \otimes k_j$$

and denote by $\mathbb{D}_{2,1}(K)$ the completion of $\nabla \varphi$ under the norm

$$\|\varphi\|_{2,1} = E^{\frac{1}{2}}\left(|\varphi|_K^2 + |\nabla \varphi|_{H \times K}^2\right).$$

For the case where $W = C_0[0,1]$, $\nabla F$ is of the form $\nabla F = \int_0^\cdot \alpha_s ds$. We will denote $\alpha_s$ by $D_s F$ and therefore for any $h = \int_0^\cdot \dot{h}_s ds$

$$(\nabla F, h)_H = \int_0^1 D_s F \cdot \dot{h}_s ds. \tag{2.3}$$

*The divergence (the Skorohod integral):* Let $v(x), x \in \mathbb{R}_n$ take values in $\mathbb{R}_n$, $v(x) = \sum_1 v_i(x)\rho_i$, where the $\rho_i$ are orthonormal vectors in $\mathbb{R}_n$. Recall, first the following "integration by parts formula":

$$\int_{\mathbb{R}_n} (v(x), \nabla F(x)) dx = -\int_{\mathbb{R}_n} F(x) \operatorname{div} v \, dx, \tag{2.4}$$

where div is the divergence: $\operatorname{div} v = \sum_1^n \frac{\partial v_i}{\partial x_i}$. Note that equation (2.4) deals with integration with respect to the Lebesgue measure on $\mathbb{R}_n$. Here we are looking for an analog of the divergence operation on $\mathbb{R}_n$ which will yield an integration by parts formula with respect to the Wiener measure.

**Definition 2.1.** *Let $u(w)$ be an $H$-valued r.v. in $(W, H, \mu)$, $u$ will be said to be in $\text{dom}_2 \delta$ if $E|u(w)|_H^2 < \infty$ and there exists a r.v. say $\delta u$ such that for all smooth functionals $f(\delta e_1, \ldots, \delta e_n)$ and all $n$ the "integration by parts" relation*

$$E\Big(\nabla f, u(w)\Big)_H = E(f \cdot \delta u) \tag{2.5}$$

*is satisfied. $\delta u$ is called the divergence or Skorohod integral.*

A necessary and sufficient condition for a square integrable $u(w)$ to be in $\text{dom}_2 \delta$ is that for some $\gamma = \gamma(u)$,

$$\left| E(u(w), \nabla f)_H \right| \leq \gamma E^{\frac{1}{2}} f^2(w)$$

for all smooth $f$. For non-random $h \in W^*$, $\delta h = \langle h, w \rangle$, setting $f = 1$ in (2.5) yields that $E\delta h = 0$. It can be shown that under proper restrictions

$$\delta\Big(f(w)u(w)\Big) = f(w)\delta u - \Big(\nabla f, u(w)\Big)_H. \tag{2.6}$$

Consequently, if $E|u|_H^2 < \infty$, and $\nabla u$ is of trace class then

$$_{W^*}\langle u, w\rangle_W = \delta u + \text{trace}\,\nabla u. \tag{2.7}$$

Among the interesting facts about the divergence operator, let us also note that for the classical Brownian motion and if $(u(w))(\cdot) = \int_0^\cdot \dot u_s(w) ds$ where $\dot u_s(w)$ is adapted and square integrable then $\delta u$ coincides with the Ito integral i.e.

$$\delta u = \int_0^1 \dot u_s(w) dw_s. \tag{2.8}$$

*Notation:* Let $u(w)$ be $H$-valued and in the domain of the divergence. We will use $\rho(u)$ to denote the exponential function

$$\rho(u) := \exp\left(\delta u - \frac{1}{2}|u|_H^2\right). \tag{2.9}$$

### 3. Realization of positive functionals by adapted shifts

Let $(W, H, \mu)$ be the classical 1-dimensional Wiener Process. The expectation operator $E$ and "a.s." will be with respect to the measure $\mu$. Let $\nu$ be a probability measure on $W$ and $\nu \ll \mu$. Set $L(w) = d\nu/d\mu$. Then by the integral representation for $L \in L^1(\mu)$

$$L(w) = 1 - \int_0^1 \dot\xi_s(w) dw_s \tag{3.1}$$

where $\dot{\xi}.$ is adapted, and $\int_0^1 \dot{\xi}_s^2(w)\,ds < \infty$ a.s. Let $L_t(w) = E(L(w)|\mathcal{F}_t)$ then

$$L_t(w) = 1 - \int_0^t \dot{\xi}_s(w)\,dw_s\,.$$

Set

$$\left.\begin{aligned}\dot{\alpha}_s(w) &= L_s^{-1}(w)\,\dot{\xi}_s(w) \quad \text{if } L_s(w) \neq 0 \\ &= 0 \quad\quad\quad\quad\quad\quad \text{if } L_s(w) = 0\end{aligned}\right\} \quad (3.2)$$

Further assume that $\int_0^1 \dot{\alpha}_s^2(w)ds < \infty$ a.s. (a sufficient condition for this is $\nu \sim \mu$). Then $L_t(w)$ solves the equation

$$L_t(w) = 1 - \int_0^t L_s(w)\,\dot{\alpha}_s(w)\,dw_s \qquad (3.3)$$

and

$$L_t(w) = \exp\left(-\int_0^t \dot{\alpha}_s(w)\,dw_s - \frac{1}{2}\int_0^t \dot{\alpha}_s^2(w)ds\right)\,. \qquad (3.4)$$

The relation between the shift $Tw = w + \int_0^\cdot \dot{\alpha}_s(w)ds$ and $L$ is given by the Girsanov theorem as

$$T(L\mu) = \mu\,. \qquad (3.5)$$

We have also the following result which assures the realization of $\nu = L \cdot \mu$:

**Proposition 3.1.** *Let $\nu, L, \dot{\alpha}$ and $T$ be as defined above. Assume that the process $w \to y(w) = (y_t(w), t \in [0,1])$ solves the equation*

$$y_t(w) = w_t - \int_0^t \dot{\alpha}_s(y)ds \qquad (3.6)$$

*in the strong sense and that*

$$\int_0^1 \dot{\alpha}_s^2(y)ds < \infty \quad\quad a.s. \qquad (3.7)$$

*Let $S : W \to W$ be the map defined by $y(w) = S(w)$, i.e., $Sw = w + \int_0^\cdot \dot{\alpha}_s(Sw)ds$. Then*

(i) *$S$ is the right inverse of $T$:*

$$T \circ S = I_W$$

*almost surely and $S$ realizes $\nu$, i.e., $S\mu = \nu$ and $S\mu \sim \mu$. In fact $T$ is almost surely invertible and $S$ is its (both sided) inverse.*

(ii) *Let $M$ be defined as*

$$M = \exp\left(\int_0^1 \dot{\alpha}_s(y)dw_s - \frac{1}{2}\int_0^1 |\dot{\alpha}_s(y)|^2 ds\right)\,,$$

then
$$\frac{dT\mu}{d\mu} = M$$

and
$$L\,M \circ T = 1 \tag{3.8}$$

$\mu$-almost surely.

*Proof.* It is clear from (3.6) that $T \circ S(w) = w$ almost surely. Moreover, it follows from the hypothesis (3.7), the random variable

$$M = \exp\left(\int_0^1 \dot{\alpha}_s(y)dw_s - \frac{1}{2}\int_0^1 |\dot{\alpha}_s(y)|^2 ds\right)$$

is almost surely finite. Using the usual stopping time technique and Fatou Lemma we have $E[M] \leq 1$. Moreover, again using the stopping technique and the Jessen theorem, cf., [23], Theorem 1.2.1, it follows that

$$E[f \circ S\,M] = E[f],$$

for any $f \in C_b(W)$. Since $M > 0$ almost surely, $S\mu$ is equivalent to $\mu$. Moreover, we have

$$L \circ S\,M = 1, \tag{3.9}$$

almost surely. To see that $S\mu = \nu$, we have, for any nice function $f$ on $W$,

$$\int_W f\,Ld\mu = \int_W f \circ S\,L \circ S\,M d\mu$$
$$= \int_W f \circ S\,d\mu,$$

by the relation (3.9). To show that $S \circ T = I_W$ almost surely, note first that $S$ is almost surely subjective: in fact, using the Girsanov theorem, we get

$$\mu(S(W)) = \int_W 1_{S(W)}\,d\mu$$
$$= \int_W 1_{S(W)} \circ S(w)\,M(w)d\mu(w)$$
$$= \int_W M\,d\mu = 1.$$

Hence $S(W) = W$ almost surely. Moreover, on the set $W_1 = \{w \in W : T \circ S(w) = w\}$ the map $S$ is injective and $\mu(W_1) = 1$ and this completes the proof of (i). For the proof of (ii), the first part follows from (i) and from the Girsanov theorem:

$$\int f \circ T\,d\mu = \int f \circ T \circ S\,Md\mu$$
$$= \int f\,Md\mu,$$

for any nice $f$ on $W$. For the second part we know already that $L \circ S M = 1$ and that $T\mu \sim \mu$ since $M > 0$ $\mu$-p.s. Hence it suffices to compose this identity with $T$. □

In the converse direction we have the following result which is based on the assumption invertibility of adapted absolutely continuous shift.

**Proposition 3.2.** *Assume that $\nu$ is realizable by*

$$Tw = w + \int_0^{\cdot} u_s(w) ds \tag{3.10}$$

*where $u.$ is adapted, $\int_0^1 u_s^2(w) ds < \infty$ a.s., $\nu = T\mu \ll \mu$ and $L = d\nu/d\mu$. Further assume that $T$ is right invertible ($TT^{-1} = I_W$ a.s.) and $\int_0^1 |u_s(T^{-1}w)|^2 ds < \infty$. Then $T$ satisfies the equation ((3.6)) with $\alpha(w) = -u(T^{-1}w)$ and*

$$L(w) = \exp\left(-\int_0^1 u_s(T^{-1}w) dw_s - \frac{1}{2}\int_0^1 |u_s(T^{-1}w)|^2 ds\right), \tag{3.11}$$

*in particular $\nu \sim \mu$.*

*Proof.* Since $T$ is right invertible, by (3.10) we have

$$T^{-1}w = w - \int_0^{\cdot} u_s(T^{-1}w) ds. \tag{3.12}$$

By Theorem 2.4.1 of [23] $(T^{-1})\mu \ll \mu$. Moreover, since $\int_0^1 |u_s|^2 ds < \infty$ a.s., it follows by Theorem 2.4.2 of [23] that $T\mu \sim \mu$. Set $\Lambda$ to be

$$\Lambda(w) = \exp\left(\int_0^1 u_s(w) dw_s - \frac{1}{2}\int_0^1 |u_s|^2(w) ds\right)$$

by the Girsanov theorem (cf.[23])

$$\Lambda = \left(\frac{dT\mu}{d\mu} \circ T\right)^{-1} = \left(L \circ T\right)^{-1}.$$

Hence $L = (\Lambda \circ T^{-1})^{-1}$ and

$$L = \exp\left(-\int_0^1 u_s \circ T^{-1} dw_s - \frac{1}{2}\int_0^1 |u_s \circ T^{-1}|^2 ds\right).$$

□

A complete characterization on the representability of Radon-Nikodym derivatives by adapted shifts is not known at present. We have, however the following result (cf. section 2.7 of [23] for the proof):

**Theorem 3.3.** *The set of measures that are realizable by adapted shifts is dense in the space of measures which are absolutely continuous with respect of the Wiener measure under the total variation norm, i.e. if on an abstract Wiener space $\nu \ll \mu$ then there exists a sequence of* invertible *and adapted shifts $T_n$ such that $dT_n\mu/d\mu$ converges to $d\nu/d\mu$ in the $L^1$ norm.*

In the rest of this section we shall make some further observations about the problem. First let us a necessary and sufficient condition for the representability of a positive random variable:

**Proposition 3.4.** *Let $U = I_W + u$, where $u : W \to H$ is adapted. Denote $\rho(-\delta u)$ (cf. equations (2.8), (2.9)) by $\Lambda$ and assume that $E[\Lambda] = 1$. Then $U\mu = L\mu$ if and only if*

$$E[\Lambda|U] \, L \circ U = 1$$

*almost surely.*

*Proof.* Necessity: from the Girsanov theorem, we have, for any $f \in C_b(W)$,

$$E[f] = E[f \circ U \, \Lambda] = E[f \circ U \, E[\Lambda|U]],$$

note that $E[\Lambda|U] = \theta \circ U$ for some measurable $\theta : W \to \mathbb{R}$. Hence we get

$$E[f] = E[f \, \theta \, L].$$

This implies that $\theta \, L = 1$ almost surely, hence $\theta \circ U \, L \circ U = E[\Lambda|U] \, L \circ U = 1$ almost surely.
Sufficiency: We have

$$\begin{aligned} E[f \circ U] &= E[f \circ U \, L \circ U \, E[\Lambda|U]] \\ &= E[f \circ U \, L \circ U \, \Lambda] \\ &= E[f \, L], \end{aligned}$$

which completes the proof. □

**Proposition 3.5.** *Assume that $L$ is given by $\rho(-\delta v)$, where $v : W \to H$ is such that $\dot{v}$ is adapted, where $v(t, w) = \int_0^t \dot{v}_s(w) ds$, $t \in [0, 1]$ and that the mapping $h \to v(w+h)$ is $\mu$-almost surely continuous and compact on the Cameron-Martin space $H$. Finally assume that*

$$\lim\sup_{|h|_H \to \infty} \frac{|v(w+h)|_H}{|h|_H} = \alpha(w),$$

*such that $0 \leq \alpha(w) < 1$ $\mu$-almost surely. Then $w \to w + v(w)$, denoted as $V(w)$, has a left inverse $U_l = I_W + u$, where $u$ is also adapted and $U_l$ realizes $L$. In fact $U_l$ is also a right inverse almost surely.*

*Proof.* The hypothesis about $v$ imply that the mapping $V$ has a left inverse $U_l$ thanks to Leray-Schauder degree theorem (cf.[23], p. 247, Theorem 9.4.1 and Remark 9.4.4). It is clear that $U_l$ is of the form $I_W + u$, with $\dot{u}$ adapted. Consequently, using the Girsanov theorem, we get, for any $f \in C_b(W)$,

$$\begin{aligned} E[f \, L] &= E[f \circ U_l \circ V \, \rho(-\delta v)] \\ &= E[f \circ U_l]. \end{aligned}$$

To show that $U_l$ is also an almost sure right inverse let us remark that it suffices to prove the subjectivity of $V$ since it is already injective almost surely by the existence of a left inverse. To show the subjectivity we need to prove that $V(W) = W$ almost surely. First the fact that $V(W)$ is measurable with respect to the universal completion of $\mathcal{B}(W)$ is given in Theorem 4.2.1 of [23]. It then follows from the Girsanov theorem that

$$\begin{aligned}\mu(V(W)) &= E[1_{V(W)} \circ V L] \\ &= E[L] = 1\,,\end{aligned}$$

since $1_{V(W)} \circ V = 1$. $\square$

We have also the following result:

**Proposition 3.6.** *Assume that $L = \rho(-\delta v)$ with $\dot{v}$ adapted and $E[L] = 1$. Let $V = I_W + v$ and represent $dV\mu/d\mu$ by $\rho(-\delta u)$ and let $U = I_W + u$. Then $U \circ V$ is a rotation, i.e., $U \circ V(\mu) = \mu$, consequently the stochastic process $(t, w) \to U \circ V(t, w)$ is a Wiener process under $\mu$, with respect to its own filtration. If it is also a martingale with respect to the canonical filtration, then $U \circ V = I_W$ almost surely. In other words $U$ is a left inverse of $V$ and it realizes $L$. In fact $U$ is also right inverse $\mu$-almost surely.*

*Proof.* For any $f \in C_b(W)$, we have

$$\begin{aligned}E[f \circ U \circ V] &= E[f \circ U\, \rho(-\delta u)] \\ &= E[f]\,,\end{aligned}$$

from the Girsanov theorem and this shows that $U \circ V$ is a rotation. Consequently it can be written as

$$U \circ V = B\,,$$

where $B$ is a Wiener process with respect to its own filtration. Besides, writing in detail the above relation, we get

$$\begin{aligned}B(w) &= V(w) + u \circ V(w) \\ &= w + v(w) + u \circ V(w)\,.\end{aligned}$$

If $B$ is also a (continuous) martingale with respect to the canonical filtration, then the adapted and finite variation process $v + u \circ V$ becomes a martingale, hence it is almost surely zero and this implies that $U \circ V = I_W$ almost surely. To complete the proof of realization of $L$ by $U$ it suffices to reason as in the proof of Proposition 3.5. The proof of the fact that $U$ is also a right inverse a.s. is also similar to that of Proposition 3.5. $\square$

**Proposition 3.7.** *Assume that $L = \rho(-\delta v)$ with $E[L] = 1$ and $\dot{v}$ adapted. Let $V = I_W + v$ assume also that $dV\mu/d\mu = \rho(-\delta u)$ and let $U = I_W + u$. If $L$ is $\sigma(V)$-measurable then $V$ is almost surely invertible and its inverse is $U$.*

*Proof.* From the Girsanov theorem, we have, for any $f \in C_b(W)$,

$$E[f \circ V \, \rho(-\delta u) \circ V \, L] = E[f \, \rho(-\delta u)] = E[f \circ V].$$

Hence we get

$$\rho(-\delta u) \circ V \, E[L|\sigma(V)] = 1$$

almost surely. The hypothesis implies that

$$\rho(-\delta u) \circ V \, L = 1$$

almost surely. Hence we have

$$\exp\left(-\delta u \circ V - \frac{1}{2}|u \circ V|_H^2\right) \exp\left(-\delta v - \frac{1}{2}|v|_H^2\right) = 1 \qquad (3.13)$$

almost surely. Note that

$$\delta u \circ V = \delta u + (u \circ V, v)_H.$$

From the equation (3.13) we get

$$\begin{aligned}
0 &= \delta u + (u \circ V, v)_H + \delta v + \frac{1}{2}(|u \circ V|_H^2 + |v|_H^2) \\
&= \delta(u+v) + \frac{1}{2}|u \circ V + v|_H^2.
\end{aligned}$$

Taking the expectation of both sides we obtain

$$u \circ V + v = 0$$

almost surely, i.e., $I_W = I_W + v + u \circ V = U \circ V$ almost surely, which means that $U$ is a left inverse of $V$. One can show that it is also a right inverse as we have already done in the preceding proofs. □

## 4. The (weak) representability of positive random variables via perturbations along the Cameron-Martin space

In this section we present the following result of Fernique [7]

**Theorem 4.1.** *Let $(W, H, \mu)$ be an abstract Wiener space. Let $\nu \ll \mu$, then there exist a probability space $(\Omega, \mathcal{H}, P)$ and three random variables $T, B$ and $Z$ defined on it such that $T, B : \Omega \to W$ and $Z : \Omega \to H$ such that $T(P) = \nu$, $B(P) = \mu$ and that $T = B + Z$.*

The results of the following lemma are not new (cf.[21, 22, 23] and the references there), we give a quick proof for the sake of completeness:

**Lemma 4.1.** *Suppose that $F: W \mapsto \mathbb{R}$ is a measurable random variable satisfying*

$$|F(w+h) - F(w)| \leq c|h|_H, \tag{4.1}$$

*almost surely, for any $h \in H$, where $c > 0$ is a fixed constant. Then $F$ belongs to $\mathbb{D}_{p,1}$ for any $p > 1$ and $\nabla F \in L^\infty(\mu, H)$. Moreover we have*

$$\mu\{|F| > \lambda\} \leq 2 \exp\left\{-\frac{(\lambda - E[F])^2}{2\|\nabla F\|^2_{L^\infty(\mu, H)}}\right\}$$

*for any $\lambda \geq 0$.*

*Proof.* We shall give the proof in the frame of a classical Wiener space, the abstract Wiener space can be reduced to this case as it is explained in Chapter 2 of [23]. Remark that there is no information about the integrability of $F$. In fact this is a consequence of the hypothesis (4.1) as explained below: Let $F_n = |F| \wedge n$, $n \in \mathbb{N}$. A simple calculation shows that

$$|F_n(w+h) - F_n(w)| \leq c|h|_H,$$

hence $F_n \in \mathbb{D}_{p,1}$ for any $p > 1$ and $|\nabla F_n| \leq c$ almost surely. We have from the Itô-Clark formula (cf.[21]),

$$F_n = E[F_n] + \int_0^1 E[D_s F_n | \mathcal{F}_s] dw_s,$$

where $D_s F$ is as defined by (3.3). From the definition of the stochastic integral, we have

$$\begin{aligned}
E\left[\left(\int_0^1 E[D_s F_n|\mathcal{F}_s]dw_s\right)^2\right] &= E\left[\int_0^1 |E[D_s F_n|\mathcal{F}_s]|^2 ds\right] \\
&\leq E\left[\int_0^1 |D_s F_n|^2 ds\right] \\
&= E[|\nabla F_n|^2] \\
&\leq c^2.
\end{aligned}$$

Since $F_n$ converges to $|F|$ in probability, and the stochastic integral is bounded in $L^2(\mu)$, by taking the difference, we see that $(E[F_n], n \in \mathbb{N})$ is a sequence of (degenerate) random variables bounded in the space of random variables under the topology of convergence in probability, denoted by $L^0(\mu)$. Therefore $\sup_n \mu\{E[F_n] > t\} \to 0$ as $t \to \infty$. Hence $\lim_n E[F_n] = E[|F|]$ is finite. Now we apply Fatou's lemma to obtain that $F \in L^2(\mu)$. Since the distributional derivative of $F$ is a bounded random variable, $F \in \mathbb{D}_{p,1}$, for any $p \geq 1$. From the above lines, replacing $F$ by $F - E[F]$, we can suppose that $E[F] = 0$. Let $(e_i) \subset H$ be a complete, orthonormal basis of $H$. Define $V_n = \sigma\{\delta e_1, \ldots, \delta e_n\}$ and let $\varphi_n = E[P_{1/n} F | V_n]$, where $P_t$ denotes the Ornstein-Uhlenbeck semigroup on $W$. Then, from Doob's Lemma,

$$\varphi_n = f_n(\delta e_1, \ldots, \delta e_n).$$

Note that, since $f_n \in \bigcap_{p,k} W_{p,k}(\mathbb{R}^n, \mu_n)$, the Sobolev embedding theorem implies that after a modification on a set of null Lebesgue measure, $f_n$ can be chosen in $C^\infty(\mathbb{R}^n)$. Let $(B_t; t \in [0,1])$ be an $\mathbb{R}^n$-valued Brownian motion. Then

$$\begin{aligned}
\mu\{|\varphi_n| > c\} &= P\{|f_n(B_1)| > c\} \\
&\leq P\{\sup_{t \in [0,1]} |E[f_n(B_1)|\mathcal{B}_t]| > c\} \\
&= P\{\sup_{t \in [0,1]} |Q_{1-t}f_n(B_t)| > c\},
\end{aligned}$$

where $P$ is the canonical Wiener measure on $C([0,1], \mathbb{R}^n)$ and $Q_t$ is the heat kernel associated to $(B_t)$, i.e.

$$Q_t(x, A) = P\{B_t + x \in A\}.$$

From the Ito formula, we have

$$Q_{1-t}f_n(B_t) = Q_1 f_n(B_0) + \int_0^t (DQ_{1-s}f_n(B_s), dB_s).$$

By definition

$$\begin{aligned}
Q_1 f_n(B_0) &= Q_1 f_n(0) = \int f_n(y) \cdot Q_1(0, dy) \\
&= \int f_n(y) e^{-\frac{|y|^2}{2}} \frac{dy}{(2\pi)^{n/2}} \\
&= E\left[E[P_{1/n}F|V_n]\right] \\
&= E\left[P_{1/n}F\right] = E[F] = 0.
\end{aligned}$$

Moreover we have $DQ_t f = Q_t Df$, hence

$$Q_{1-t}f_n(B_t) = \int_0^t (Q_{1-s}Df_n(B_s), dB_s) = M_t^n.$$

The Doob-Meyer process $(\langle M^n, M^n \rangle_t, t \in \mathbb{R}_+)$ of the martingale $M^n$ can be controlled as

$$\begin{aligned}
\langle M^n, M^n \rangle_t &= \int_0^t |DQ_{1-s}f_n(B_s)|^2 ds \\
&\leq \int_0^t \|Df_n\|_{C_b}^2 ds = t\|\nabla f_n\|_{C_b}^2 \\
&= t\|\nabla f_n\|_{L^\infty(\mu_n)} \\
&\leq t\|\nabla F\|_{L^\infty(\mu, H)}^2.
\end{aligned}$$

Hence from the exponential Doob inequality, we obtain

$$P\left\{\sup_{t \in [0,1]} |Q_{1-t}f_n(B_t)| > \lambda\right\} \leq 2\exp\left\{-\frac{\lambda^2}{2\|\nabla F\|_{L^\infty(\mu, H)}^2}\right\}.$$

Consequently

$$\mu\{|\varphi_n| > \lambda\} \leq 2\exp\left\{-\frac{\lambda^2}{2\|\nabla F\|^2_{L^\infty(\mu,H)}}\right\}.$$

Since $\varphi_n \to \varphi$ in probability the proof is completed. □

The following proposition will play the intermediate role in proof of Theorem 4.1:

**Proposition 4.2.** *Let $\mu$ be the standard Gaussian measure on $\mathbb{R}^n$ and let $\nu$ be another probability on $\mathbb{R}^n$ such that $\nu \ll \mu$ and that the Radon-Nikodym derivative $\frac{d\nu}{d\mu}$ is essentially bounded by a constant $A < \infty$. Then there exist a Gaussian random variable $B$ with values in $\mathbb{R}^n$, and an $\mathbb{R}^n$-valued random variable $Z$, both defined on $(\Omega, \mathcal{H}, P)$, satisfying*

$$E[|Z|] \leq \sqrt{\frac{\pi}{2}} A$$

*and the probability measure on $\mathbb{R}^n$ induced by $T = B + Z$ is $\nu$.*

*Proof.* For the proof we shall use basically the Kantorovitch-Rubinstein theorem: Let $g : \mathbb{R}^n \to \mathbb{R}$ be Lipschitz with Lipschitz constant 1. By Lemma 4.1, and denoting by $E_\mu, E_\nu$ the expectation with respect to the $\mu$ and $\nu$ measures respectively:

$$\mu\left\{\left|g - E_\mu[g]\right| > t\right\} \leq 2\exp-\frac{t^2}{2}.$$

Therefore

$$\nu\left\{\left|g - E_\mu[g]\right| > t\right\} \leq A\exp-\frac{t^2}{2}. \tag{4.2}$$

Note that if for some r.v. $Y$ it holds that $\psi(t) = \text{Prob}\{|Y| > t\} < Ae^{\frac{-t^2}{2}}$ then, from the Fubini theorem

$$E[|Y|] = \int_0^\infty \psi(t)dt \leq A\int_0^\infty e^{-\frac{t^2}{2}}dt = A\sqrt{\frac{\pi}{2}}.$$

Therefore by (4.2)

$$\left|E_\nu[g] - E_\mu[g]\right| \leq E_\nu\left|g - E_\mu[g]\right| \leq A\sqrt{\frac{\pi}{2}}. \tag{4.3}$$

Define, for any two measures $\alpha$ and $\beta$ on $W$

$$K(\alpha, \beta) = \sup_{\|g\|_{\text{Lip}} \leq 1}\left|\int g d\alpha - \int g d\beta\right|$$

where the supremum is over all Lipschitz functions with Lipschitz constant 1. From (4.3)

$$K(\mu, \nu) \leq A\sqrt{\frac{\pi}{2}}. \tag{4.4}$$

Now let $N(\mu, \nu)$ be the Rubinstein distance between $\mu$ and $\nu$:

$$N(\mu, \nu) = \inf \int_{W \times W} |x - y|_W \, M(dx, dy) \tag{4.5}$$

where the infimum is taken over the set $\Sigma$ which consists of all the probability measures $M$ on $W \times W$ for which $M(dx, W) = \mu(dx)$ and $M(W, dy) = \nu(dy)$. Note that the set $\Sigma$ is compact w.r.to the weak topology of the measures on $W \times W$ and the integral at the r.h.s. of (4.5) is lower semicontinuous under this topology. Hence there is at least one measure $M_0$ which realizes the infimum provided that the infimum is finite. Moreover, by the Kantorovitch-Rubinstein theorem (cf. e.g. [6, 11.8.2]), $K(\mu, \nu) = N(\mu, \nu)$ and this implies the finiteness of $N(\mu, \nu)$ hence also the existence of $M_0(dx, dy)$ on $W \times W$ with the prescribed properties. Letting $\Omega = W \times W$, $P = M_0$, $B = \Pi_1$, $T = \Pi_2$ and $Z = T - B$, where $\Pi_i$, $i = 1, 2$ are the projection maps from $W \times W$ to $W$, we see that $T(P) = \nu$, $B(P) = \mu$ and that $T = B + Z$. □

**Proof of Theorem 4.1:** assume first that $L = \frac{d\nu}{d\mu} \leq A < \infty$. Let $(e_i, i \geq 1) \subset W^*$ be a complete orthonormal basis of $H$, set

$$\pi_N(w) = w_N = \sum_1^N \langle w, e_i \rangle e_i$$

and let $\mu_N, \nu_N$ denote the measures on $\mathbb{R}^N$ which are the images of $\mu$ and $\nu$ under $\pi_N$. Let $T_N$ and $B_N$ be as in the proof of the proposition, then $(T_N, N \geq 1)$ and $(B_N, N \geq 1)$ converge in law and by (4.4)

$$E\left[\left|T_N - B_N\right|_H\right] \leq A\sqrt{\frac{\pi}{2}}.$$

It follows by weak convergence and the Skorohod representation theorem (cf. e.g. theorem 7.1.4. and sections 11.5, 11.7 of [6]), there exists a probability space $(\Omega, \mathcal{H}, P)$ and that there exists a sequence of random variables $(T'_N, B'_N, N \in \mathbb{N})$ which converges a.s. to $(T', B')$ such that the laws of $T'$ and $B'$ under the probability $P$ are $\nu$ and $\mu$ respectively and moreover, by Fatou's lemma and by the lower semi continuity of the Cameron-Martin norm with respect to the topology of $W$,

$$E_P\left[\left|T' - B'\right|_H\right] \leq A\sqrt{\frac{\pi}{2}}.$$

In order to complete the proof it remains to remove the restriction about the boundedness of the density $L$. Let $A_0 = 0$, $A_n = 2^{n-1}$ and $W_n = \{w : L(w) \in [A_{n-1}, A_n)\}$. Set $L_n(w) = 1_{W_n} L(w)$. Set

$$d\nu_n = \begin{cases} \frac{L_n d\mu}{\int L_n d\mu} & \text{if } \int L_n d\mu > 0 \\ d\mu & \text{otherwise} \end{cases}.$$

Applying the results for $\mu, \nu_n$ as above yields, for any $n \geq 1$, $(T_n'', B_n'')$ defined on $(\Omega_n, \mathcal{H}_n, P_n)$ such that $T_n''(P_n) = \nu_n$ and $B_n''(P_n) = \mu$ with

$$E_{P_n}\left[\left|T_n'' - B_n''\right|_H\right] \leq \left(\frac{A_n}{\int L_n d\mu}\right)\sqrt{\frac{\pi}{2}}$$

if $\int L_n d\mu > 0$. Let $\Omega$ be the disjoint union of $(\Omega_n, n \geq 1)$ and define

$$P = \sum_n \left(\int L_n d\mu\right) \times P_n.$$

Define, for $\omega \in \Omega_n$

$$T''(\omega) = T_n''(\omega), \quad B''(\omega) = B_n''(\omega)$$

then $T''(P) = \nu$ and $B''(P) = \mu$ and $Z'' = T'' - B'' \in H$ almost surely. $\square$

## 5. Realization by triangular transformations

### 5.1. Triangular transformations

A map $T : \mathbb{R}^n \to \mathbb{R}^n$ is said to be *triangular* if

$$T(x_1, \ldots, x_n) = \begin{pmatrix} T_1(x_1) \\ T_2(x_1, x_2) \\ \vdots \\ T_n(x_1, \ldots, x_n) \end{pmatrix}$$

where $T_i : \mathbb{R}^i \to \mathbb{R}^1$. A triangular map is called *increasing* if for every $i = 1, \ldots, n$, $x_i \to T_i(x_1, \ldots, x_i)$ is non decreasing.

**Proposition 5.1.** *Let $\mu$ and $\nu$ be probability measures on $\mathbb{R}^n$. Assume that $\mu$ and $\nu$ are absolutely continuous with respect to the Lebesgue measure on $\mathbb{R}^n$, then there exists an increasing triangular transformation $T$ on $\mathbb{R}^n$ such that*

$$T\mu = \nu.$$

*Proof.* We proceed by induction on dimension: (a) Consider first the case $n = 1$. Let $F_\mu(x)$ and $F_\nu(x)$, $x \in \mathbb{R}^1$ denote the distribution functions of $\mu$ and $\nu$ respectively. Assume, that $F_\nu$ is strictly increasing, then for every $x \in \mathbb{R}^1$ there is a unique $y \in \mathbb{R}^1$ such that $F_\mu(x) = F_\nu(y)$. Set $T_{\mu,\nu}(x) = y$ then $F_\nu(T_{\mu,\nu}x) = F_\mu(x)$. Since $F_\nu$ is strictly increasing, $F_\nu$ is invertible and

$$T_{\mu,\nu}x = F_\nu^{-1}(F_\mu x) \tag{5.1}$$

transforms $\mu$ into $\nu$ and $T_{\mu,\nu}$ is strictly increasing. If $F_\nu$ is not strictly increasing, set

$$F_\nu^{-1}(y) := \inf\left(\beta : F_\nu(\beta) \geq y\right).$$

Note that with this interpretation of $F_\nu^{-1}$, $T_{\mu,\nu}$ is increasing and equation (5.1) proves the proposition for $n = 1$.

(b) Consider now the case $n = 2$, let $\varphi(x_1, x_2)$ and $\psi(x_1, x_2)$ denote the probability density with respect to the Lebesgue measure of $\mu$ and $\nu$ respectively. Let $\varphi^{x_1}(x_2)$ and $\psi^{x_1}(x_2)$ denote the regular probability densities of $x_2$ conditioned on $x_1$, i.e.

$$\varphi(x_1, x_2) = \varphi_1(x_1)\varphi^{x_1}(x_2)$$

and

$$\psi(x_1, x_2) = \psi_1(x_1)\psi^{x_1}(x_2).$$

Let $T_1 : \mathbb{R}^1 \to \mathbb{R}^1$ denote the transformation on $\mathbb{R}^1$ transforming $d\mu_1 = \varphi_1(x_1)dx_1$ to the measure $d\nu_1 = \psi(y_1)dy_1$ as constructed in part (a). Next, let $T_2(x_1, x_2)$ be the map transforming, for fixed $x_1$, the measure $\varphi^{x_1}(x_2)dx_2$ to the measure $\psi^{x_1}(x_2)dx_2$. Then

$$T = \begin{pmatrix} T_1 \\ T_2 \end{pmatrix}$$

is an increasing triangular mapping and $T\mu = \nu$.

(c) For $n > 2$ the result follows by induction using the conditional density of $x_n$ given $x_1, \ldots, x_{n-1}$. □

It follows from the proof by induction that

**Corollary 5.2.** *Proposition 5.1 holds for $\mathbb{R}^\infty$, provided that the projections of the measures $\mu$ and $\nu$ to $\mathbb{R}^n$ are absolutely continuous with respect to the Lebesgue measure for any $n \geq 1$.*

### *5.2. The realization of a Radon-Nikodym derivative by a triangular transformation:*

Let $(W, H, \mu)$ be an AWS. Let $(e_i, i \geq 1) \subset W^*$ be an orthonormal basis[2] for $H$. By the Ito-Nisio theorem (cf. e.g. [15] or [23]) $\left\| w - \sum_{i=1}^N {}_{W^*}\langle e_i, w\rangle_W \right\|_W \xrightarrow{a.s.} 0$ as $N \to \infty$. Let $\nu \ll \mu$ then as shown in the previous section we can construct an increasing triangular transformation $T$:

$$Tw = \sum_i T_i(w)e_i$$

$$T_iw = g_i\left({}_{W^*}\langle e_1, w\rangle_W, \ldots, {}_{W^*}\langle e_i, w\rangle_W\right) \qquad (5.2)$$

such that the measure induced by $Tw$ on $W$ is $\nu$. This yields a measurable transformation which induces $\nu$ and yields another proof to the result of Section 4 that even if $\nu \ll \mu$ then $\nu$ is realizable.

---
[2]We do not distinguish $e_i \in W^*$ from its image in $H$ if there is no confusion.

### *5.3. The entropy associated with a triangular transformation:*

Let $\mu$ denote the standard Gaussian measure on $\mathbb{R}^n$ and let $T$ denote a $C^1$ diffeomorphic map on $\mathbb{R}^n$. By the change of variables formula, it holds that for every $g \in C_b(\mathbb{R}^n)$

$$\int_{\mathbb{R}^n} g \circ T(x) |\Lambda(x)| d\mu(x) = \int_{\mathbb{R}^n} g(x) d\mu(x) \tag{5.3}$$

where $\Lambda$ is given as follows: set $T(x) = x + u(x)$ then

$$\Lambda(x) = \det(I + \nabla u) \exp\left\{-(u(x), x)_{\mathbb{R}^n} - \frac{1}{2}|u(x)|^2_{\mathbb{R}^n}\right\} \tag{5.4}$$

(cf. e.g. [23]). The Carleman-Fredholm determinant of an $n \times n$ matrix $I_{\mathbb{R}^n} + A$, $\det_2(I_{\mathbb{R}^n} + A)$, is defined as

$$\det_2(I_{\mathbb{R}^n} + A) = \det(I_{\mathbb{R}^n} + A) e^{-\operatorname{trace} A}$$

and $\Lambda$ can be rewritten as

$$\Lambda = \det_2(I + \nabla u) \exp\left\{\operatorname{trace} \nabla u - (u(x), x) - \frac{1}{2}|u(x)|^2_{\mathbb{R}^n}\right\} \tag{5.5}$$

and by (2.8):

$$\Lambda = \det_2(I + \nabla u) \exp\left\{\delta u - \frac{1}{2}|u(x)|^2\right\}. \tag{5.6}$$

The reason of rewriting $\Lambda$ with the modified Carleman-Fredholm[3] determinant lies in the fact that $\det_2(I+A)$ is well-defined for the Hilbert-Schmidt operators, in fact $A \to \det_2(I + A)$ is analytic on the space of Hilbert-Schmidt operators on an infinite dimensional Hilbert space, while $A \to \det(I + A)$ is defined only for the nuclear (or trace class) operators. Assume now that the measure induced by $T$ is absolutely continuous with respect to $\mu$, $T\mu \ll \mu$, then

$$\int_{\mathbb{R}^n} g \circ T(x) d\mu(x) = \int_{\mathbb{R}^n} f(x) L(x) d\mu(x) \tag{5.7}$$

where $L$ is the Radon-Nikodym derivative. Hence comparing (5.7) with (5.3) and recalling that $T$ is *invertible* we have

$$L(x) = \left|\Lambda(T^{-1}x)\right|^{-1} \tag{5.8}$$

where $\Lambda$ is as defined by (5.3).

Let $L$ be the Radon-Nikodym derivative of $\nu$ with respect to $\mu$ and assume that $L \log L \in L^1(\mu)$ set $L(x) \log L(x) = 0$ whenever $L(x) = 0$. The relative entropy $\mathcal{H}(\nu|\mu)$, of $\nu$ with respect to $\mu$ is defined as

$$\mathcal{H}(\nu|\mu) = E_\mu(L \log L) = E_\nu(\log L).$$

---

[3]cf. Chapter 3 of [23] for details.

**Proposition 5.3.** *Let $T$ be a $C^1$ and invertible, triangular and increasing map on $\mathbb{R}^n$, let $\mu$ be the standard Gaussian measure on $\mathbb{R}^n$, assume that $T\mu \ll \mu$ and $L = dT\mu/d\mu$ satisfies $L \log L \in L^1(\mu)$. Then, setting $T(x) = x + u(x)$, it holds that*

$$\frac{1}{2} E_\nu[|u|^2] \leq \mathcal{H}(\nu|\mu). \tag{5.9}$$

*Proof.* By (5.8) and by (5.7) and then (5.5)

$$\begin{aligned}
\int_{\mathbb{R}^n} L(x) \log L(x) d\mu(x) &= \int_{\mathbb{R}^n} L(x) \log \left|\Lambda(T^{-1}x)\right|^{-1} d\mu(x) \\
&= \int_{\mathbb{R}^n} -\log|\Lambda(x)| d\mu(x) \\
&= -\int_{\mathbb{R}^n} \left\{ \log \left|\det{}_2(I_{\mathbb{R}^n} + \nabla u(x))\right| \right. \\
&\qquad \left. + \operatorname{trace} \nabla u(x) - (u(x), x) - \frac{1}{2}|u(x)|^2 \right\} d\mu(x).
\end{aligned}$$

Therefore we get the simple form for $\mathcal{H}(\nu|\mu)$:

$$\begin{aligned}
E_\mu[L \log L] &= \frac{1}{2} E_\mu[|u|^2] - E_\mu\left[\log |\det{}_2(I + \nabla u)|\right] \\
&= \frac{1}{2} E_\mu[|u|^2] + E_\mu\left[\log \det{}_2(I + \nabla u)\right]
\end{aligned}$$

where the last line follows since $\det_2(I + \nabla u) \in [0, 1]$ almost surely as explained below: $Tx = x + u(x)$ is triangular and increasing, therefore $\nabla u$ is a triangular matrix whose diagonal elements are all non-negative. Hence setting $\frac{du_i(x)}{dx_i} = \alpha_i \geq 0$, yields

$$\det{}_2(I + \nabla u) = \prod_{i=1}^n (1 + \alpha_i) e^{-\alpha_i}.$$

Since for $\alpha \geq 0$, $0 \leq (1 + \alpha) e^{-\alpha} \leq 1$ it follows that $0 \leq \det_2(I + \nabla u) \leq 1$, and (5.9) follows. □

### 5.4. The realization by $H$-valued shifts

**Theorem 5.4.** *Let $(W, H, \mu)$ be an abstract Wiener space, define a probability measure $\nu$ on $(W, \mathcal{B}(W))$ as $d\nu = L d\mu$ where $L \geq 0$ with $E[L] = 1$. Then there exist a measurable $T$ on $W$ such that $T\mu = \nu$ and $\mu$-almost surely $Tw - w \in H$, i.e. $T$ is a perturbation of identity.*

*Proof.* Assume first that $L \log L \in L^1(\mu)$. Let $(e_i, i \geq 1)$ be as in subsection 5.2, let $\mu_n$ and $\nu_n$ be the images of probability measures $\mu$ and $\nu$ induced on $\mathbb{R}^n$ by the map $w \to \pi_n(w) = (\langle e_1, w \rangle, \ldots, \langle e_n, w \rangle)$. Then $\nu_n \ll \mu_n$, set $l_n = d\nu_n/d\mu_n$ then $l_n \circ \pi_n = L_n$ is the conditional expectation of $L$ conditioned on the $\sigma$-field generated by $w \to \pi_n(w) = (\langle e_1, w\rangle, \ldots, \langle e_n, w\rangle)$. Since $s \to s\log s$ is a convex function on $\mathbb{R}_+$, it follows by Jensen's inequality that

$$E_\mu[L_n \log L_n] \leq E_\mu[L \log L] < \infty.$$

Let $t_n$ be the triangular map such that $t_n\mu_n = \nu_n$, define $T_n = t_n \circ \pi_n$, then from the construction of subsection 5.2, $\lim_n T_n = T$ almost surely. Setting $u_n(w) = T_n(w) - w$, we have from Proposition 5.3 that for all $n$

$$\sup_n E_\mu\left[|u_n|_H^2\right] \leq 2E_\mu[L \log L]$$

hence by Fatou's lemma $E[|Tw - w|_H^2] < \infty$.

In order to remove the restriction $L \log L \in L_1(\mu)$ we proceed as follows: partition the space $W$ into sets $E_m = \{w : m \leq L(w) < m+1\}$, $m \in \mathbb{N}$, avoiding sets on which $\nu(E_m) = 0$. Set $\nu^m(A) = \nu(A \cap E_m)$ for all Borel sets $A$ in $W$. Next partition the real line by intervals $D_m$ such that $\mu_1(D_m) = \nu(E_m)$ where $\mu_1$ is the standard 1-dimensional Gaussian measure on $\mathbb{R}e_1$. Consider the decomposition $\mu = \mu_1 \otimes \mu'$ where $\mu'$ is the (Gaussian) measure defined as $(I_W - \pi_1)(\mu)$. Define $\mu_{1,m}$ on the real line as $\mu_{1,m}(B) = \mu(B \cap D_m)$ for $B \in \mathcal{B}(\mathbb{R})$. Then

$$\mu = \sum_m \mu_{1,m} \otimes \mu'.$$

Define $t_m : \mathbb{R}e_1 \to \mathbb{R}e_1$ such that $t_m(\mu_{1,m}) = \nu(E_1)\mu_1$. Then $(t_m \times I_{W'})(\mu_{1,m} \otimes \mu') = \nu(E_m)\mu$. Transform now $\nu(E_m)\mu$ in to $\nu^m$ as in the first part of this subsection with a triangular transformation, say $U_m$ and define

$$T_m = U_m \circ (t_m \times I_{W'})$$

where $W' = \pi_1^{-1}\{0\}$. Define finally $T$ as to be $T_m$ on the set $D_m \times W'$. Note that $T_m(w) - w$ is a.s. $H$-valued. Since for $m_1 \neq m_2$, the domains of $T_{m_1}$ and $T_{m_2}$ are disjoint and so are their corresponding ranges since they are contained in the supports of $\nu^{m_1}$ and $\nu^{m_2}$ respectively. Consequently $(T_m, m \geq 1)$ yields a triangular transformation transforming $\mu$ to $\nu$ and such that $Tw - w$ in a.s. in the Cameron-Martin space. □

## 6. The realization by gradient shifts

Let $W$ be a separable Fréchet space with its Borel sigma algebra $\mathcal{B}(W)$ and assume that there is a separable Hilbert space $H$ which is injected densely and continuously into $W$, thus the topology of $H$ is, in general, stronger than the topology induced by $W$. The cost function $c : W \times W \to \mathbb{R}_+ \cup \{\infty\}$ is defined as

$$c(x, y) = |x - y|_H^2,$$

we suppose that $c(x, y) = \infty$ if $x - y$ does not belong to $H$. Clearly, this choice of the function $c$ is not arbitrary, in fact it is closely related to Ito Calculus, hence also to the problems originating from Physics, quantum chemistry, large deviations, etc. Since for all the interesting measures on $W$, the Cameron-Martin space is a negligible set, the cost function will be infinity very frequently. Let $\Sigma(\rho,\nu)$ denote the set of probability measures on $W \times W$ with given marginals $\rho$ and $\nu$. It is a convex, compact set under the weak topology $\sigma(\Sigma, C_b(W \times W))$. The problem of Monge consists of finding a measurable map $T: W \to W$, called the optimal transport of $\rho$ to $\nu$, i.e., $T\rho = \nu$ which minimizes the total cost

$$U \to \int_W |x - U(x)|_H^2 d\rho(x),$$

between all the maps $U: W \to W$ such that $U\rho = \nu$. On the other hand the Monge-Kantorovitch problem consists of finding a measure on $W \times W$, which minimizes the function $\theta \to J(\theta)$, defined by

$$J(\theta) = \int_{W \times W} |x - y|_H^2 d\theta(x,y), \tag{6.1}$$

where $\theta$ runs in $\Sigma(\rho,\nu)$. Note that $\inf\{J(\theta) : \theta \in \Sigma(\rho,\nu)\}$ is the square of Wasserstein metric $d_H(\rho,\nu)$ with respect to the Cameron-Martin space $H$. Any solution $\gamma$ of the Monge-Kantorovitch problem will give a solution to the Monge problem provided that its support is included in the graph of a map.

Let us recall some notions of convexity on the Wiener space (cf. [8, 23]). Let $K$ be a Hilbert space, a subset $S$ of $K \times K$ is called cyclically monotone if any finite subset $\{(x_1, y_1), \ldots, (x_N, y_N)\}$ of $S$ satisfies the following algebraic condition:

$$\langle y_1, x_2 - x_1 \rangle + \langle y_2, x_3 - x_2 \rangle + \cdots + \langle y_{N-1}, x_N - x_{N-1} \rangle + \langle y_N, x_1 - x_N \rangle \leq 0,$$

where $\langle \cdot, \cdot \rangle$ denotes the inner product of $K$. It turns out that $S$ is cyclically monotone if and only if

$$\sum_{i=1}^{N} \langle y_i, x_{\sigma(i)} - x_i \rangle \leq 0,$$

for any permutation $\sigma$ of $\{1, \ldots, N\}$ and for any finite subset $\{(x_i, y_i) : i = 1, \ldots, N\}$ of $S$. Note that $S$ is cyclically monotone if and only if any translate of it is cyclically monotone. By a theorem of Rockafellar, any cyclically monotone set is contained in the graph of the subdifferential of a convex function in the sense of convex analysis ([19]) and even if the function may not be unique its subdifferential is unique.

A measurable function defined on $(W, H, \mu)$ with values in $\mathbb{R} \cup \{\infty\}$ is called 1-convex if the map

$$h \to f(x+h) + \frac{1}{2}|h|_H^2 = F(x,h)$$

is convex on the Cameron-Martin space $H$ with values in $L^0(\mu)$. Note that this notion is compatible with the $\mu$-equivalence classes of random variables thanks to the Cameron-Martin theorem. It is proven in [8] that this definition is equivalent the following condition: Let $(\pi_n, n \geq 1)$ be a sequence of regular, finite dimensional, orthogonal projections of $H$, increasing to the identity map $I_H$. Denote also by $\pi_n$ its continuous extension to $W$ and define $\pi_n^\perp = I_W - \pi_n$. For $x \in W$, let $x_n = \pi_n x$ and $x_n^\perp = \pi_n^\perp x$. Then $f$ is 1-convex if and only if

$$x_n \to \frac{1}{2}|x_n|_H^2 + f(x_n + x_n^\perp)$$

is $\pi_n^\perp \mu$-almost surely convex.

**Definition 6.1.** *Let $\xi$ and $\eta$ be two probabilities on $(W, \mathcal{B}(W))$. We say that a probability $\gamma$ on $(W \times W, \mathcal{B}(W \times W))$ is a solution of the Monge-Kantorovitch problem associated to the couple $(\xi, \eta)$ if the first marginal of $\gamma$ is $\xi$, the second one is $\eta$ and if*

$$J(\gamma) = \int_{W \times W} |x - y|_H^2 d\gamma(x, y) = \inf \left\{ \int_{W \times W} |x - y|_H^2 d\beta(x, y) : \beta \in \Sigma(\xi, \eta) \right\},$$

*where $\Sigma(\xi, \eta)$ denotes the set of all the probability measures on $W \times W$ whose first and second marginals are respectively $\xi$ and $\eta$. We shall denote the Wasserstein distance between $\xi$ and $\eta$, which is the positive square-root of this infimum, with $d_H(\xi, \eta)$.*

**Remark:** By the weak compactness of probability measures on $W \times W$ and the lower semi-continuity of the strictly convex cost function, the infimum in the definition is attained even if the functional $J$ is identically infinity.

The following result, whose proof is outlined below (cf. also[9, 10]) is an extension of an inequality due to Talagrand [20] and it gives a sufficient condition for the Wasserstein distance to be finite:

**Theorem 6.1.** *Suppose that $L$ is a positive random variable with $L \log L \in L^1(\mu)$ and with $E[L] = 1$. Let $\nu$ be the measure $d\nu = Ld\mu$, we then have*

$$d_H^2(\nu, \mu) \leq 2E[L \log L]. \tag{6.2}$$

*Proof.* Let us remark first that we can take $W$ as the classical Wiener space $W = C_0([0, 1])$ and, using the stopping techniques of the martingale theory, we may assume that $L$ is upper and lower bounded almost surely. Then a classical result of the Ito calculus implies that $L$ can be represented as an exponential martingale

$$L_t = \exp\left\{-\int_0^t \dot{u}_\tau dW_\tau - \frac{1}{2} \int_0^t |\dot{u}_\tau|^2 d\tau\right\},$$

with $L = L_1$. Let us define $u : W \to H$ as $u(t, x) = \int_0^t \dot{u}_\tau d\tau$ and $U : W \to W$ as $U(x) = x + u(x)$. The Girsanov theorem implies that $x \to U(x)$ is a Brownian motion under $\nu$, hence $\beta = (U \times I)\nu \in \Sigma(\mu, \nu)$. Let $\gamma$ be any optimal measure, then

$$\begin{aligned} J(\gamma) &= d_H^2(\nu, \mu) \leq \int_{W \times W} |x - y|_H^2 d\beta(x, y) \\ &= E[|u|_H^2 L] \\ &= 2E[L \log L], \end{aligned}$$

where the last equality follows also from the Girsanov theorem and the Ito stochastic calculus. $\square$

Combining Theorem 6.1 with the triangle inequality for the Wasserstein distance gives:

**Corollary 6.2.** *Assume that $\nu_i$ ($i = 1, 2$) have Radon-Nikodym densities $L_i$ ($i = 1, 2$) with respect to the Wiener measure $\mu$ which satisfy $E|L_i(w) \log L_i(w)| < \infty, i = 1, 2$. Then*

$$d_H(\nu_1, \nu_2) < \infty.$$

Let us give a simple application of the above result

**Corollary 6.3.** *Assume that $A \in \mathcal{B}(W)$ is any set of positive Wiener measure. Define the $H$-gauge function of $A$ as*

$$q_A(x) = \inf(|h|_H : h \in (A - x) \cap H).$$

*Then we have*

$$E[q_A^2] \leq 2 \log \frac{1}{\mu(A)},$$

*in other words*

$$\mu(A) \leq \exp\left\{-\frac{E[q_A^2]}{2}\right\}.$$

*Similarly if $A$ and $B$ are $H$-separated, i.e., if $A_\varepsilon \cap B = \emptyset$, for some $\varepsilon > 0$, where $A_\varepsilon = \{x \in W : q_A(x) \leq \varepsilon\}$, then*

$$\mu(A_\varepsilon^c) \leq \frac{1}{\mu(A)} e^{-\varepsilon^2/4}$$

*and consequently*

$$\mu(A)\, \mu(B) \leq \exp\left(-\frac{\varepsilon^2}{4}\right).$$

**Remark 6.1.** We already know that, from the $0-1$–law, $q_A$ is almost surely finite, besides it satisfies $|q_A(x + h) - q_A(x)| \leq |h|_H$, hence the hypothesis of Lemma 4.1 are satisfied. Consequently $E[\exp \lambda q_A^2] < \infty$ for any $\lambda < 1/2$ (cf. the Appendix B.8 of [23] and [22]). In fact all these assertions can also be proved with the technique used below.

*Proof.* Let $\nu_A$ be the measure defined by

$$d\nu_A = \frac{1}{\mu(A)} 1_A d\mu \,.$$

Let $\gamma_A$ be the solution of the Monge-Kantorovitch problem, it is easy to see that the support of $\gamma_A$ is included in $W \times A$, hence

$$|x - y|_H \geq \inf\{|x - z|_H : z \in A\} = q_A(x)\,,$$

$\gamma_A$-almost surely. This implies in particular that $q_A$ is almost surely finite. It follows now from the inequality (6.2)

$$E[q_A^2] \leq -2 \log \mu(A)\,,$$

hence the proof of the first inequality follows. For the second let $B = A_\varepsilon^c$ and let $\gamma_{AB}$ be the solution of the Monge-Kantorovitch problem corresponding to $\nu_A, \nu_B$. Then we have from the Corollary 6.2,

$$d_H^2(\nu_A, \nu_B) \leq -4 \log \mu(A)\mu(B)\,.$$

Besides the support of the measure $\gamma_{AB}$ is in $A \times B$, hence $\gamma_{AB}$-almost surely $|x - y|_H \geq \varepsilon$ and the proof follows. □

We call optimal every probability measure[4] $\gamma$ on $W \times W$ such that $J(\gamma) < \infty$ and that $J(\gamma) \leq J(\theta)$ for every other probability $\theta$ having the same marginals as those of $\gamma$. We recall that a finite dimensional subspace $F$ of $W$ is called regular if the corresponding projection is continuous. Similarly a finite dimensional projection of $H$ is called regular if it has a continuous extension to $W$.

The proof of the next theorem, for which we refer the reader to [10], can be done by choosing a proper disintegration of any optimal measure in such a way that the elements of this disintegration are the solutions of finite dimensional Monge-Kantorovitch problems. The latter is proven with the help of the section-selection theorem [4, 5].

The notion of spread measure is key to all the results about the measure transportation:

**Definition 6.2.** *A probability measure $m$ on $(W, \mathcal{B}(W))$ is called a spread measure if there exists a sequence of finite dimensional regular projections $(\pi_n, n \geq 1)$ converging to $I_H$ such that the regular conditional probabilities $m(\cdot | \pi_n^\perp = x_n^\perp)$, (where $\pi_n^\perp = I_W - \pi_n$) which are concentrated in the n-dimensional spaces $\pi_n(W) + x_n^\perp$ vanish on the sets of Hausdorff dimension $n - 1$ for $\pi_n^\perp(m)$-almost all $x_n^\perp$ and for any $n \geq 1$.*

**Remark 6.2.** *Clearly any measure absolutely continuous with respect to $\mu$ is spread.*

---

[4]In fact the results of this section are true for the positive measures of equal mass.

**Theorem 6.4 (General Monge-Kantorovitch transportation).** *Suppose that $\rho$ and $\nu$ are two probability measures on $W$ such that*

$$d_H(\rho, \nu) < \infty$$

*and that $\rho$ is spread. Then there exists a unique solution of the Monge-Kantorovitch problem, denoted by $\gamma \in \Sigma(\rho, \nu)$ and $\gamma$ is supported by the graph of a Borel map $T$ which is the solution of the Monge problem. $T : W \to W$ is of the form $T = I_W + \xi$, where $\xi \in H$ almost surely. Besides we have*

$$\begin{aligned} d_H^2(\rho, \nu) &= \int_{W \times W} |T(x) - x|_H^2 d\gamma(x, y) \\ &= \int_W |T(x) - x|_H^2 d\rho(x), \end{aligned}$$

*and, with the notations of Definition 6.2, for $\pi_n^\perp \rho$-almost almost all $x_n^\perp$, the map $u \to \xi(u + x_n^\perp)$ is cyclically monotone on $(\pi_n^\perp)^{-1}\{x_n^\perp\}$, in the sense that*

$$\sum_{i=1}^N \left(\xi(x_n^\perp + u_i), u_{i+1} - u_i\right)_H \leq 0$$

$\pi_n^\perp \rho$-*almost surely, for any cyclic sequence $\{u_1, \ldots, u_N, u_{N+1} = u_1\}$ from $\pi_n(W)$. Finally, if $\nu$ is also a spread measure then $T$ is invertible, i.e, there exists $S : W \to W$ of the form $S = I_W + \eta$ such that $\eta \in H$ satisfies a similar cyclic monotononicity property as $\xi$ and that*

$$\begin{aligned} 1 &= \gamma\{(x, y) \in W \times W : T \circ S(y) = y\} \\ &= \gamma\{(x, y) \in W \times W : S \circ T(x) = x\}. \end{aligned}$$

*In particular we have*

$$\begin{aligned} d_H^2(\rho, \nu) &= \int_{W \times W} |S(y) - y|_H^2 d\gamma(x, y) \\ &= \int_W |S(y) - y|_H^2 d\nu(y). \end{aligned}$$

The case where one of the measures is the Wiener measure and the other is absolutely continuous with respect to $\mu$ is the most important one for the applications. Consequently we give the related results separately in the following theorem where the tools of the Malliavin calculus give more information about the maps $\xi$ and $\eta$ of Theorem 6.4:

**Theorem 6.5 (Gaussian case).** *Let $\nu$ be the measure $d\nu = Ld\mu$, where $L$ is a positive random variable, with $E[L] = 1$. Assume that $d_H(\mu, \nu) < \infty$ (for instance $L \in \mathbb{L} \log \mathbb{L}$). Then there exists a 1-convex function $\phi \in \mathbb{D}_{2,1}$, unique up to a constant, such that the map $T = I_W + \nabla \phi$ is the unique solution of the* *original problem of Monge. Moreover, its graph supports the unique solution of the Monge-Kantorovitch problem $\gamma$. Consequently*

$$(I_W \times T)\mu = \gamma$$

*In particular $T$ maps $\mu$ to $\nu$ and $T$ is almost surely invertible, i.e., there exists some $T^{-1}$ such that $T^{-1}\nu = \mu$ and that*

$$\begin{aligned} 1 &= \mu\{x : T^{-1} \circ T(x) = x\} \\ &= \nu\{y \in W : T \circ T^{-1}(y) = y\}\,. \end{aligned}$$

**Remark 6.3.** *Assume that the operator $\nabla$ is closable with respect to $\nu$, then we have $\eta = \nabla\psi$. In particular, if $\nu$ and $\mu$ are equivalent, then we have*

$$T^{-1} = I_W + \nabla\psi\,,$$

*where is $\psi$ is a 1-convex function.*

Let us give some applications of the above theorem to the factorization of the absolutely continuous transformations of the Wiener measure.

Assume that $V = I_W + v : W \to W$ be an absolutely continuous transformation and let $L \in \mathbb{L}^1_+(\mu)$ be the Radon-Nikodym derivative of $V\mu$ with respect to $\mu$. Let $T = I_W + \nabla\phi$ be the transport map such that $T\mu = L.\mu$. Then it is easy to see that the map $s = T^{-1} \circ V$ is a rotation, i.e., $s\mu = \mu$ (cf. [23]) and it can be represented as $s = I_W + \alpha$. In particular we have

$$\alpha + \nabla\phi \circ s = v\,. \tag{6.3}$$

Since $\phi$ is a 1-convex map, we have $h \to \frac{1}{2}|h|_H^2 + \phi(x+h)$ is almost surely convex (cf.[8]). Let $s' = I_W + \alpha'$ be another rotation with $\alpha' : W \to H$. By the 1-convexity of $\phi$, we have

$$\frac{1}{2}|\alpha'|_H^2 + \phi \circ s' \geq \frac{1}{2}|\alpha|_H^2 + \phi \circ s + \left(\alpha + \nabla\phi \circ s, \alpha' - \alpha\right)_H,$$

$\mu$-almost surely. Taking the expectation of both sides, using the fact that $s$ and $s'$ preserve the Wiener measure $\mu$ and the identity (6.3), we obtain

$$E\left[\frac{1}{2}|\alpha|_H^2 - (v,\alpha)_H\right] \leq E\left[\frac{1}{2}|\alpha'|_H^2 - (v,\alpha')_H\right]\,.$$

Hence we have proven the existence part of the following

**Proposition 6.6.** *Let $\mathcal{R}_2$ denote the subset of $L^2(\mu, H)$ whose elements are defined by the property that $x \to x + \eta(x)$ is a rotation, i.e., it preserves the Wiener measure. Then $\alpha$ is the unique element of $\mathcal{R}_2$ which minimizes the functional*

$$\eta \to M_v(\eta) = E\left[\frac{1}{2}|\eta|_H^2 - (v,\eta)_H\right]\,.$$

*Proof.* The only claim to prove is the uniqueness and it follows easily from Theorem 6.5. □

The following theorem, whose proof is rather easy, gives a better understanding of the structure of absolutely continuous transformations of the Wiener measure:

**Theorem 6.7.** *Assume that $U : W \to W$ be a measurable map and $L \in \mathbb{L}\log\mathbb{L}(\mu)$ a positive random variable with $E[L] = 1$. Assume that the measure $\nu = L \cdot \mu$ is a Girsanov measure for $U$, i.e., that one has*

$$E[f \circ U\, L] = E[f],$$

*for any $f \in C_b(W)$. Then there exists a unique map $T = I_W + \nabla\phi$ with $\phi \in \mathbb{D}_{2,1}$ is 1-convex, and a measure preserving transformation $R : W \to W$ such that $U \circ T = R$ $\mu$-almost surely and $U = R \circ T^{-1}$ $\nu$-almost surely.*

Another version of Theorem 6.7 can be announced as follows:

**Theorem 6.8.** *Assume that $Z : W \to W$ is a measurable map such that $Z\mu \ll \mu$, with $d_H(Z\mu, \mu) < \infty$. Then $Z$ can be decomposed as*

$$Z = T \circ s,$$

*where $T$ is the unique transport map of the Monge-Kantorovitch problem for $\Sigma(\mu, Z\mu)$ and $s$ is a rotation.*

Here is a more general version of the polar factorization:

**Theorem 6.9.** *Assume that $\rho$ and $\nu$ are spread measures with $d_H(\rho, \nu) < \infty$ and that $U\rho = \nu$, for some measurable map $U : W \to W$. Let $T$ be the optimal transport map sending $\rho$ to $\nu$, whose existence and uniqueness is proven in Theorem 6.4. Then $R = T^{-1} \circ U$ is a $\rho$-rotation (i.e., $R\rho = \rho$) and $U = T \circ R$, moreover, if $U$ is a perturbation of identity, then $R$ is also a perturbation of identity. In both cases, $R$ is the $\rho$-almost everywhere unique minimal $\rho$-rotation in the sense that*

$$\int_W |U(x) - R(x)|_H^2 d\rho(x) = \inf_{R' \in \mathcal{R}} \int_W |U(x) - R'(x)|_H^2 d\rho(x), \qquad (6.4)$$

*where $\mathcal{R}$ denotes the set of $\rho$-rotations.*

*Proof.* Let $T$ be the optimal transportation of $\rho$ to $\nu$ whose existence and uniqueness follows from Theorem 6.4. The unique solution $\gamma$ of the Monge-Kantorovitch problem for $\Sigma(\rho, \nu)$ can be written as $\gamma = (I \times T)\rho$. Since $\nu$ is spread, $T$ is invertible on the support of $\nu$ and we have also $\gamma = (T^{-1} \times I)\nu$. In particular $R\rho = T^{-1} \circ U\rho = T^{-1}\nu = \rho$, hence $R$ is a rotation. Let $R'$ be another rotation in $\mathcal{R}$, define $\gamma' = (R' \times U)\rho$, then $\gamma' \in \Sigma(\rho, \nu)$ and the optimality of $\gamma$ implies that $J(\gamma) \leq J(\gamma')$, besides we have

$$\begin{aligned}\int_W |U(x) - R(x)|_H^2 d\rho(x) &= \int_W |U(x) - T^{-1} \circ U(x)|_H^2 d\rho(x) \\ &= \int_W |x - T^{-1}(x)|_H^2 d\nu(x) \\ &= \int_W |T(x) - x|_H^2 d\rho(x) \\ &= J(\gamma).\end{aligned}$$

On the other hand

$$J(\gamma') = \int_W |U(x) - R'(x)|_H^2 d\rho(x),$$

hence the relation (6.4) follows. Assume now that for the second rotation $R' \in \mathcal{R}$ we have the equality

$$\int_W |U(x) - R(x)|_H^2 d\rho(x) = \int_W |U(x) - R'(x)|_H^2 d\rho(x).$$

Then we have $J(\gamma) = J(\gamma')$, where $\gamma'$ is defined above. By the uniqueness of the solution of Monge-Kantorovitch problem due to Theorem 6.4, we should have $\gamma = \gamma'$. Hence $(R \times U)\rho = (R' \times U)\rho = \gamma$, consequently, we have

$$\int_W f(R(x), U(x))d\rho(x) = \int_W f(R'(x), U(x))d\rho(x),$$

for any bounded, measurable map $f$ on $W \times W$. This implies in particular

$$\int_W (a \circ T \circ R)(b \circ U)d\rho = \int_W (a \circ T \circ R')(b \circ U)d\rho$$

for any bounded measurable functions $a$ and $b$. Let $U' = T \circ R'$, then the above expression reads as

$$\int_W a \circ U \, b \circ U d\rho = \int_W a \circ U' \, b \circ U d\rho.$$

Taking $a = b$, we obtain

$$\int_W (a \circ U)(a \circ U')\, d\rho = \|a \circ U\|_{L^2(\rho)} \|a \circ U'\|_{L^2(\rho)},$$

for any bounded, measurable $a$. This implies that $a \circ U = a \circ U'$ $\rho$-almost surely for any $a$, hence $U = U'$ i.e., $T \circ R = T \circ R'$ $\rho$-almost surely. Let us denote by $S$ the left inverse of $T$ whose existence follows from Theorem 6.4 and let $D = \{x \in W : S \circ T(x) = x\}$. Since $\rho(D) = 1$ and since $R$ and $R'$ are $\rho$-rotations, we have also

$$\rho\left(D \cap R^{-1}(D) \cap R'^{-1}(D)\right) = 1.$$

Let $x \in W$ be any element of $D \cap R^{-1}(D) \cap R'^{-1}(D)$, then

$$\begin{aligned} R(x) &= S \circ T \circ R(x) \\ &= S \circ T \circ R'(x) \\ &= R'(x), \end{aligned}$$

consequently $R = R'$ on a set of full $\rho$-measure. □

Although the following result is a translation of the results of this section, it is interesting from the point of view of stochastic differential equations:

**Theorem 6.10.** *Let $(W, \mu, H)$ be the standard Wiener space on $\mathbb{R}^d$, i.e., $W = C(\mathbb{R}_+, \mathbb{R}^d)$. Assume that there exists a probability $P \ll \mu$ which is the weak solution of the stochastic differential equation*

$$dy_t = dW_t + b(t, y)dt,$$

*such that $d_H(P, \mu) < \infty$. Then there exists a process $(T_t, t \in \mathbb{R}_+)$ which is a pathwise solution of some (anticipative) stochastic differential equation whose law is equal to $P$.*

*Proof.* Let $T$ be the transport map constructed in Theorem 6.5 corresponding to $dP/d\mu$. Then it has an inverse $T^{-1}$ such that $\mu\{T^{-1} \circ T(x) = x\} = 1$. Let $\phi$ be the 1-convex function such that $T = I_W + \nabla\phi$ and denote by $(D_s\phi, s \in \mathbb{R}_+)$ the representation of $\nabla\phi$ in $L^2(\mathbb{R}_+, ds)$. Define $T_t(x)$ as the trajectory $T(x)$ evaluated at $t \in \mathbb{R}_+$. Then it is easy to see that $(T_t, t \in \mathbb{R}_+)$ satisfies the stochastic differential equation

$$T_t(x) = W_t(x) + \int_0^t l(s, T(x))ds, \ t \in \mathbb{R}_+,$$

where $W_t(x) = x(t)$ and $l(s, x) = D_s\phi \circ T^{-1}(x)$ if $x \in T(W)$ and zero otherwise. □

**Definition 6.3.** *Assume that $\rho$ and $\nu$ are two probability measures on $(W, \mathcal{B}(W))$. We say that they are at locally finite Wasserstein distance from each other if there exists a measurable partition $(Z_n, n \geq 1)$ of $W$ such that*

$$d_H(\rho, \nu_n) < \infty,$$

*for any $n \geq 1$, where $\nu_n$ is the probability defined by $\nu_n(A) = \nu(A \cap Z_n)/\nu(Z_n)$.*

**Remark 6.4.** *Any measure $\nu$ which is absolutely continuous with respect to $\mu$ is at locally finite distance from $\mu$.*

**Theorem 6.11.** *Assume that $\rho$ and $\nu$ are at locally finite distance from each other and that $\rho$ is a spread measure. Then there exists a map $T : W \to W$, of the form $T = I_W + \xi$ (i.e., perturbation of identity) such that $T$ maps the measure $\rho$ to the measure $\nu$.*

*Proof.* Let $(Z_n, n \geq 1)$ be a partition of $W$ such that $d_H(\rho, \nu_n) < \infty$ for each $n \in \mathbb{N}$. Let $(M_n, n \geq 1)$ be a partition of $W$ such that $\rho(M_n) = \nu(Z_n)$ for each $n \in \mathbb{N}$. Define $t_n$ as the map which maps the measure $1_{M_n} \rho$ to the measure $\nu(Z_n) \rho$. Let also $U_n$ be the $H$-cyclically monotone map whose existence has been proven in [10], which maps $\nu(Z_n)\rho$ to the measure $1_{Z_n} \nu$. Define $T$ on $M_n$ as the map $T_n = U_n \circ t_n$. Then it is easy to see that $T\rho = \nu$:

$$\begin{aligned}
\int f \circ T d\rho &= \sum_n \int_{M_n} f \circ U_n \circ t_n \, d\rho \\
&= \sum_n \nu(Z_n) \int f \circ U_n \, d\rho \\
&= \sum_n \int_{Z_n} f \, d\nu \\
&= \int f \, d\nu .
\end{aligned}$$

□

**Question:** We do not know in general that given two measures one of which being spread, if there is just one perturbation of identity which is cyclically monotone which maps the spread one to the second measure. Even in the case where the first measure is $\mu$ and the second one is $d\nu = Ld\mu$ where $L$ is positive and $E[L] = 1$ this question is unanswered. In finite dimensional case the answer to this question is positive; in the infinite dimensional case the difficulty is due to the fact that the Cameron-Martin norm on $W$ is only lower semi-continuous, hence, the weak limit of measures of cyclically monotone support on $W \times W$ do not have necessarily cyclically monotone support.

### 7. The realization by a "signal and independent noise"

#### *7.1. Introduction*

Let $(W, H, \mu)$ be a standard one-dimensional Wiener space and let $(\Theta, \mathcal{F}, \rho)$ be a probability space unrelated to $(W, H, \mu)$. Let $(W \times \Theta, \mathcal{B}(W) \vee \mathcal{F}, \mu \times \rho)$ be the corresponding product space. Let $(u(\theta), \theta \in \Theta)$ be a random variable on $\Theta$ taking values in the Cameron-Martin space $H$, i.e. $u = \int_0^{\cdot} \dot{u}_s(\theta) \, ds$. Let $y = u(\theta) + w$ and as specified above $w$ and $u$ are independent. We will use $\mu_{Y|U}$ to denote the conditional measure of $y$ conditioned on the sub sigma field generated by $u$; also, $\mu_{U|Y}$ will denote the conditional measure of $u$ conditioned on the subsigma field generated by $y$. Then

$$\frac{d\mu_{Y|U}}{d\mu}(w, u) = \exp\left(\int_0^1 \dot{u}_s(\theta) dw_s - \frac{1}{2} \int_0^1 |\dot{u}_s(\theta)|^2 ds\right)$$

and
$$L(w) = \frac{d\mu_Y}{d\mu}(w) = \int_H \frac{d\mu_{Y|U}}{d\mu}(w, u(\theta))\rho(d\theta) \qquad (7.1)$$

The problem considered in this section is the characterization of those non-negative random variables on $(W, H, \mu)$, say $L(w)$, where $E[L] = 1$ which are realizable by the model of equation (7.1). We will assume throughout this section that $E[L^2] < \infty$.

### 7.2. Preliminaries

Since $E[L^2] < \infty$, we have by the Wiener-Ito representation
$$L(w) = 1 + \sum_{n=1}^{\infty} \delta^n(f_n)$$

where $(f_n(t_1, \cdots, t_n), n \geq 1)$, are real valued and symmetric functions on $[0,1]^n$, $\delta^n(f^n)$ is the multiple Wiener-Ito integral
$$\int_{[0,1]^n} f(t_1, \cdots, t_n)\, dw(t_1) \ldots dw(t_n)$$

and
$$E[L^2] = 1 + \sum_1^{\infty} n!\|f_n\|_{L^2([0,1]^n)}^2.$$

A random variable $F(w)$ will be said to belong to $\mathcal{H}_\lambda$ for some real $\lambda \geq 1$ if $E[F^2] < \infty$ and $E[F_\lambda^2] < \infty$ where
$$F_\lambda(w) = \sum_1^{\infty} \lambda^n\, \delta^n(f_n).$$

For $0 < \lambda < 1$, $\mathcal{H}_\lambda$ is defined as the dual to $\mathcal{H}_{1/\lambda}$ (relative to $L^2(\mu)$) and in this case the elements of $\mathcal{H}_\lambda$ are generalized Wiener functionals. Set
$$\mathcal{H}_0 = \bigcup_{\lambda > 0} \mathcal{H}_\lambda\,;\ \mathcal{H}_\infty = \bigcap_{\lambda > 1} \mathcal{H}_\lambda.$$

Then, for $\lambda > 1$
$$\mathcal{H}_\infty \subset \mathcal{H}_\lambda \subset \mathcal{H}_1 = L^2(W) \subset \mathcal{H}_{1/\lambda} \subset \mathcal{H}_0.$$

An element $F \in \mathcal{H}_0$ will be said to be *positive* and denoted $F \geq 0$ if $\langle F, G \rangle \geq 0$ for every element $G$ of $\mathcal{H}_\infty$ and satisfying $G \geq 0$. For $F \in \mathcal{H}_1$ this means that $F(w)$ is $\mu$-a.s. non negative. Moreover, let $P_\alpha, \alpha > 0$, denote the Ornstein-Uhlenbeck semigroup
$$P_\alpha F(w) = E\left[F\left(e^{-\alpha}w + \sqrt{1 - e^{-2\alpha}}(\tilde{w})\right)|\mathcal{B}(W)\right]$$

where $\tilde{w}$ is an independent copy of the Brownian motion $w$ living on $(\tilde{W}, \mathcal{B}(\tilde{W}))$. Therefore if $F = \sum \delta^n(f_n)$, then

$$P_\alpha F = \sum e^{-n\alpha} \delta^n(f_n)$$

and since $P_\alpha$ preserves positivity, if $F \in \mathcal{H}_0$ and $F_\lambda \geq 0$ then $F_\beta \geq 0$ for all $\beta$ satisfying $0 < \beta < \lambda$. In fact $P_\alpha$ is even positivity improving in the sense that if $F \in L^1(\mu)$ is a.s. non-negative then $P_\alpha F > 0$ almost surely.

### 7.3. The positivity of random variables

Let $F$ be in $\mathcal{H}_0$, define $\lambda_{\text{pos}}(F)$, *the positivity index* of $F$, as the supremum of all $\lambda > 0$ such that $F_\lambda \geq 0$ and $\lambda_{\text{pos}}(F) = 0$ if this set is empty. Hence $F_\lambda \geq 0$ for all $\lambda \in (0, \lambda_{\text{pos}}(F))$.

In order to continue we recall the following two definitions:

(a) Let $\phi(y)$ be a complex-valued on $y \in H$. $\phi(y)$ is said to be *positive definite* if for all $n$, all $y_i \in H$, $C_i \in \mathbb{C}$, $i = 1, 2, \ldots n$

$$\sum_{i,j=1}^n C_i C_j^* \phi(y_i - y_j) \geq 0.$$

(b) Let $F \in \mathcal{H}_0$ and for any $h \in H$ set $\varepsilon(h) = \exp(i\delta h + \frac{1}{2}|h|_H^2)$. The modified $\tau$ transform of $F$, $\tilde{\tau}_F$, is defined by

$$\tilde{\tau}_F(h) = \langle F, \varepsilon(h) \rangle = \sum_{n=0}^\infty i^n \langle f_n, h^{\otimes n} \rangle_{H^{\otimes n}}. \tag{7.2}$$

It follows directly from the last equation that

$$\tilde{\tau}_{F_\lambda}(h) = \tilde{\tau}_F(\lambda h).$$

The following is a key result for the characterization of the realization problem in the next subsection.

**Theorem 7.1.** *Let $F \in \mathcal{H}_0$, then $F \geq 0$ if and only if $h \to \tilde{\tau}_F(h) \exp -\frac{1}{2}|h|_H^2$ is positive definite on $H$.*

*Proof.* The necessity is straightforward: Suppose that $F \geq 0$. Then for any $\{h_1, \cdots, h_n\} \subset H$, $\{c_1, \cdots, c_n\} \subset \mathbb{C}$ and $n \geq 1$, we have

$$\sum_{i,j=1}^n c_i \bar{c}_j \tau_f(h_i - h_j) e^{-\frac{1}{2}\|h_i - h_j\|^2} = \sum_{i,j=1}^n c_i \bar{c}_j E\left[F e^{\sqrt{-1}\delta(h_i - h_j)}\right]$$

$$= E\left[F \left|\sum_{j=1}^n e^{\sqrt{-1}\delta h_j} c_j\right|^2\right], \tag{7.3}$$

where $G = \left|\sum_{j=1}^n e^{\sqrt{-1}\delta h_j} c_j\right|^2 \geq 0$ and $G \in \mathcal{H}_\infty = \cap_\lambda \mathcal{H}_\lambda$; therefore $\langle F, G \rangle \geq 0$ and hence $h \to \tilde{\tau}_F(h) \exp -\frac{1}{2}\|h\|^2$ is positive definite. For the converse direction cf. [16] or chapter 5 of [13] and the references therein. □

**Remark:** It is easy to verify from the last theorem that $F \in \mathcal{H}_0$ is *strongly positive* if and only if $\tilde{\tau}_F(h)$ is positive definite.

### 7.4. The realization by a "signal and independent noise" model

Let $(\Theta, \mathcal{F}, \rho)$ be as defined earlier, let $u$ be an $H$-valued random variable on this space. Further assume that for all $n \geq 1$ $t_1, \ldots, t_n \subset [0, 1]$, $E[u(t_1)u(t_2)\cdots u(t_n)]$ exists as a symmetric square integrable symmetric function $f_n(t_1, \cdots, t_n) \in L^2([0,1]^n)$. The collection of all such $u$ defined on some probability spaces $(\Theta, \mathcal{F}_\theta, \rho)$ will be denoted by $U_H$.

**Definition 7.1.** *A generalized functional $F \in \mathcal{H}_0$ will be said to be* induced by $u$, *if there exists a probability space $(\Theta, \mathcal{F}, \rho)$ and an $H$ valued r.v. $u$ defined on it such that*
$$f_n(t_1, \cdots, t_n) = \frac{1}{n!} E[u(t_1)u(t_2)\cdots u(t_n)]$$
*in other words*
$$f_n = \frac{1}{n!} E\left[u^{\otimes n}\right]$$
*for all $n \geq 1$. Now, if $F \in L^2(\mu)$ is induced by some $u \in U_H$ and*
$$\sum_{n=1}^\infty (n!)^{-1} E\left[\|u^{\otimes n}\|_{H^{\otimes n}}^2\right] < \infty$$
*then, setting*
$$Y(w, \theta) = \exp\left(\delta u - \frac{1}{2}|u|_H^2\right)$$
*yields*
$$F(w) = E\left[\exp \delta u - \frac{1}{2}|u|_H^2 \Big| \mathcal{F}_W\right]$$
$$= \sum_{n=0}^\infty E\left[\delta^n u^{\otimes n} | \mathcal{F}_W\right]$$
$$= \sum_{n=0}^\infty \delta^n \left(E_\rho[u^{\otimes n}]\right)$$

We can now state the following characterization of the realizability of Radon-Nikodym derivatives by a "signal plus independent noise" model.

**Theorem 7.2.** *Let $F \in \mathcal{H}_0$ with $F = \sum_{n=0}^{\infty} \delta^n(f_n)$ and $\langle F, 1 \rangle = 1$ then there exists a cylindrical probability measure $p$ on $H$ such that for all $h_1, \cdots, h_n \in H$, $n \geq 1$*

$$\langle f_n, h_1 \otimes \cdots \otimes h_n \rangle_{H^{\otimes n}} = \frac{1}{n!} \int_H (h, h_1)_H \ldots (h, h_n)_H \, p(dh)$$

*and $f_2(t_1, t_2)$ in of trace class on $[0,1]^2$ if and only if $F$ is realizable.*

*Proof.* If a probability measure exists, then cf. equation (7.2) and theorem 7.1

$$\int_H \exp i(h_0, h)_H p(dh) = 1 + \sum_1^{\infty} i^n \langle f_n, h^{\otimes n} \rangle = \tilde{\tau}(h).$$

Therefore $\tilde{\tau}_F$ is positive definite hence $F$ is strongly positive. The converse follows by Theorem 7.1 and Theorem 2.1 of [12] since the trace-class property of $f_2$ implies the fact that the measure $p$ is in fact sigma additive measure on $H$. This means the existence of an $H$-valued random variable $u$ whose law is $p$. In fact, let $(\Theta, \mathcal{F}, \rho)$ be any probability space rich enough to induce a random variable $u : \Theta \to H$ whose law is $p$. Define on the product space $(W \times \Theta, \mathcal{B}(W) \otimes \mathcal{F}, \mu \times \rho)$ the map $Y = u(\theta) + w$, then clearly

$$\frac{d\mu_Y}{d\mu}(w) = F(w).$$

□